\documentclass[11pt, reqno]{amsart}

\usepackage{latexsym}
\usepackage{amssymb}
\usepackage{mathrsfs}
\usepackage{amsmath}
\usepackage{fancybox,color}
\usepackage{enumerate}
\usepackage[latin1]{inputenc}

\newcommand{\kom}[1]{}
 \def\1{\raisebox{2pt}{\rm{$\chi$}}}

\newtheorem{theorem}{Theorem}[section]

\newtheorem{remark}[theorem]{Remark}

\newcommand{\R}{{\mathbb R}}

 \newcommand{\eps}{{\varepsilon}}
 \def\1{\raisebox{2pt}{\rm{$\chi$}}}
 
\newcommand{\abs}[1]{\left|#1\right|}

\newcommand{\Rn}{\mathbb{R}^n}

\def\vint_#1{\mathchoice%
          {\mathop{\kern 0.2em\vrule width 0.6em height 0.69678ex depth -0.58065ex
                  \kern -0.8em \intop}\nolimits_{\kern -0.4em#1}}%
          {\mathop{\kern 0.1em\vrule width 0.5em height 0.69678ex depth -0.60387ex
                  \kern -0.6em \intop}\nolimits_{#1}}%
          {\mathop{\kern 0.1em\vrule width 0.5em height 0.69678ex
              depth -0.60387ex
                  \kern -0.6em \intop}\nolimits_{#1}}%
          {\mathop{\kern 0.1em\vrule width 0.5em height 0.69678ex depth -0.60387ex
                  \kern -0.6em \intop}\nolimits_{#1}}}
\def\vintslides_#1{\mathchoice%
          {\mathop{\kern 0.1em\vrule width 0.5em height 0.697ex depth -0.581ex
                  \kern -0.6em \intop}\nolimits_{\kern -0.4em#1}}%
          {\mathop{\kern 0.1em\vrule width 0.3em height 0.697ex depth -0.604ex
                  \kern -0.4em \intop}\nolimits_{#1}}%
          {\mathop{\kern 0.1em\vrule width 0.3em height 0.697ex depth -0.604ex
                  \kern -0.4em \intop}\nolimits_{#1}}%
          {\mathop{\kern 0.1em\vrule width 0.3em height 0.697ex depth -0.604ex
                  \kern -0.4em \intop}\nolimits_{#1}}}

\newcommand{\kint}{\vint}
\newcommand{\intav}{\vint}
\newcommand{\aveint}[2]{\mathchoice%
          {\mathop{\kern 0.2em\vrule width 0.6em height 0.69678ex depth -0.58065ex
                  \kern -0.8em \intop}\nolimits_{\kern -0.45em#1}^{#2}}%
          {\mathop{\kern 0.1em\vrule width 0.5em height 0.69678ex depth -0.60387ex
                  \kern -0.6em \intop}\nolimits_{#1}^{#2}}%
          {\mathop{\kern 0.1em\vrule width 0.5em height 0.69678ex depth -0.60387ex
                  \kern -0.6em \intop}\nolimits_{#1}^{#2}}%
          {\mathop{\kern 0.1em\vrule width 0.5em height 0.69678ex depth -0.60387ex
                  \kern -0.6em \intop}\nolimits_{#1}^{#2}}}
\newcommand{\ud}{\, d}
\newcommand{\half}{{\frac{1}{2}}}

\newcommand{\ol}{\overline}
\newcommand{\Om}{\Omega}
\newcommand{\I}{\textrm{I}}
\newcommand{\II}{\textrm{II}}

\newcommand{\om}{\omega}

\newcommand{\spt}{\operatorname{spt}}

\begin{document}

\title[Regularity for nonlinear stochastic games]
{Regularity for nonlinear stochastic games}

\author[Luiro]{Hannes Luiro}
\address{Department of Mathematics and Statistics, University of
Jyv\"askyl\"a, PO~Box~35, FI-40014 Jyv\"askyl\"a, Finland,
hannes.s.luiro@jyu.fi}
\email{hannes.s.luiro@jyu.fi}

\author[Parviainen]{Mikko Parviainen}
\address{Department of Mathematics and Statistics, University of
Jyv\"askyl\"a, PO~Box~35, FI-40014 Jyv\"askyl\"a, Finland,
mikko.j.parviainen@jyu.fi}
\email{mikko.j.parviainen@jyu.fi}

\thanks{Both authors have been supported by the Academy of Finland.} 
\subjclass[2010]{91A15, 35J92, 35B65, 35J60, 49N60}
\keywords{Dynamic programming principle, $p$-Laplace, tug-of-war, tug-of-war with noise with space dependent probabilities}
\begin{abstract}
We establish regularity for functions satisfying a dynamic programming equation, which may arise for example from stochastic games or  discretization schemes.  Our results can also be utilized in obtaining regularity and existence results for the corresponding partial differential equations.  
\end{abstract}

\maketitle

\section{Introduction}


In \cite{luirops13}, we studied  regularity for the stochastic game called tug-of-war with noise. Recalling the recently discovered connection to the  $p$-harmonic functions \cite{peress08}, our results implied local Lipschitz regularity  for the solutions to the $p$-Laplace equation
for $2< p<\infty$.    The approach was based on a choice of
strategies  for the players, and is thus quite different from the PDE proofs.

Our argument utilized symmetry properties of strategies, and a sharp cancellation effect produced by this symmetry, which directly allowed us  to obtain  a local Lipschitz estimate. It is a nontrivial task to extend this method to a more general class of problems where the perfect symmetry breaks down. Thus, in this paper, we develop a more flexible regularity method. As a starting point, we take a dynamic programming equation 
\[
\begin{split}
u(x) =\,\frac{1}{2}\sup_{\mu_1\in\mathcal{A}_1(x)}\int_{\Rn}u(y)\,d\mu_1(y)\,
+\frac{1}{2}\inf_{\mu_2\in\mathcal{A}_2(x)}\int_{\Rn}u(y)\,d\mu_2(y)\,
\end{split}
\]
as explained in detail in Section \ref{sec:applicability}.
This is a rather general formulation that covers a wide class of games; in addition to stochastic games it may as well arise from discretization schemes in numerical methods for partial differential equations, see for example \cite{oberman05}. 

To illustrate our approach, we prove (asymptotic) Hölder continuity, Theorem \ref{thm:main},  for several examples using the method. In particular, we prove the regularity for the tug-of-war with noise with space dependent probabilities, Section \ref{sec:p(x)}. In addition to the break down of the symmetry, the problem in applying known approaches is the lack of translation invariance. Furthermore, we study the tug-of-war with noise related to the $p$-Laplacian with the full range $1<p\le \infty$ in Section \ref{sec:1<p}. 
Proving a local regularity result for this game may seem complicated, since the players can affect the direction of the noise, and thus it is not easy to say much about the noise distribution. However, the method of this paper is well suited for the task.   We did not exhaust the list of possible examples that can be treated by the method but expect it to be useful in many more problems. Also, at least in some cases, the method can be improved to give directly stronger regularity results. 

Our method  arises from stochastic game theory even if for expository reasons we have eventually avoided stochastic arguments. The idea is that we start the game simultaneously at two points $x\in \Rn$ and $z\in \Rn$, and try to pull the points 'closer' to each other. Here  closer means, at least roughly, in the sense of  averages and in terms of a suitable comparison function. 
 To show that we may pull the points closer in this sense, we may consider the process in the higher dimensional space by setting  $(x,z)\in \R^{2n}$, and use the subspace
$$
T:=\{(x,z)\in \R^{2n}\,:\, x=z\}\subset \R^{2n}
$$ 
as a target. We use the following strategy: if our opponent takes a non optimal step we pull directly towards $T$. If the opponent pulls (almost) away from $T$, then our strategy is to aim at the exactly opposite step. The curvature of the comparison function gives an advantage to us.
It is also worth noting that there is a freedom to choose among the probability measures in $\R^{2n}$ having the measures arising from the original games as marginals; cf.\ the setting in the optimal mass transport problems. Suitable choices will be helpful in the proofs.

After finishing the paper it has come to our attention that couplings of stochastic processes have been employed in the study of regularity for second order linear uniformly parabolic equations with continuous highest order coefficients, see for example \cite{lindvallr93}, \cite{priolaw06}, and \cite{kusuoka15}. The method here has also some similarities to the Ishii-Lions method \cite{ishiil90}, see also for example \cite{porrettap13}. However, we do not rely on the theorem of sums in the theory of viscosity solution, but the proofs are built on the ideas arising from the game theory. 

Our motivation to study the above problems is threefold: First, the study of stochastic games has received a lot of attention  on their own right because of deep mathematical questions that arise and also due to their central role in many applications. Second, the dynamic programming principle can be interpreted as a discretization of the PDE. Thus, results can also be interpreted as results for the corresponding numerical schemes. Third, by passing to the limit with respect to the step size our results imply regularity results for the PDE. 



There is a powerful connection between  the classical linear partial differential equations and probability theory. 
In the nonlinear case, the connection between the games and Bellman-Isaacs equations was established in the 80s. However, a similar connection between the normalized
$p$-Laplace or $\infty$-Laplace equations (which are discontinuous operators in the gradient variable) and the tug-of-war games with noise was discovered only rather recently in \cite{peress08,peresssw09}. This connection has later been extended or utilized in several different contexts, see for example  \cite{atarb10, manfredipr10,bjorklandcf12a, bjorklandcf12, manfredipr12,antunovicpss12, lius15, ruosteenoja16}. 

\section{The regularity method}\label{prelim}

\subsection{Background}
\label{sec:applicability}

Let  $\Omega\subset\Rn,\ n\ge 2$, be an open and bounded set. Denote by $\mathcal{M}$ the space of unit Radon measures on $\Rn$, and suppose that we are given
\[
\begin{split}
\mathcal{A}:=(\mathcal A_1,\mathcal A_2)&:\Rn\to \mathcal{P}(\mathcal{M})\times\mathcal{P}(\mathcal{M})\,,
\end{split}
\]
where $\mathcal{P}(\mathcal{M})$ denotes the family of all subsets of $\mathcal{M}\,$. 
Our objective is to develop a regularity method for  bounded Borel measurable solutions $u:\Omega \to \R$ to the dynamic programming equation
\begin{equation}
\label{DPP}
u(x) =\,\frac{1}{2}\sup_{\mu_1\in\mathcal{A}_1(x)}\int_{\Rn}u(y)\,d\mu_1(y)\,
+\frac{1}{2}\inf_{\mu_2\in\mathcal{A}_2(x)}\int_{\Rn}u(y)\,d\mu_2(y)\,.
\end{equation}
The heuristic idea is that $u$ represents a value function for a two players zero-sum game, where the rules of the game are determined by $\mathcal{A}_1(x)$ and $\mathcal{A}_2(x)$. More precisely, when the game token is at $x\in\Omega$, the maximizer can choose $\mu_1\in\mathcal{A}_1(x)$ and minimizer $\mu_2\in\mathcal{A}_2(x)$ and then the distribution of the next location is determined by the probability measure $\frac{\mu_1+\mu_2}{2}\,$.  
Then we can write
\begin{align*}
u(x)-u(z)&=\frac{1}{2}\sup_{\mu_1\in\mathcal{A}_1(x)}\int_{\Rn}u(y)\,d\mu_1(y)\,
+\frac{1}{2}\inf_{\mu'_2\in\mathcal{A}_2(x)}\int_{\Rn}u(y)\,d\mu'_2(y)\,
\\
&\hspace{1 em}-\frac{1}{2}\sup_{\mu'_1\in\mathcal{A}_1(z)}\int_{\Rn}u(y)\,d\mu'_1(y)\,
-\frac{1}{2}\inf_{\mu'_2\in\mathcal{A}_2(z)}\int_{\Rn}u(y)\,d\mu'_2(y)\,\\
&=\sup_{\mu_1\in\mathcal{A}_1(x),\mu_2\in\mathcal{A}_2(z)}\bigg[\inf_{\mu'_2\in\mathcal{A}_2(x),\mu'_1\in\mathcal{A}_1(z)}\\
&\hspace{3 em}\bigg(\int_{\Rn}u(y)\,d\big(\frac{\mu_1+\mu'_2}{2}\big)(y)\,-\int_{\Rn}u(y)\,d\big(\frac{\mu'_1+\mu_2}{2}\big)(y)\bigg)\bigg]\,.
\end{align*}
Let $\Gamma(\frac{\mu_1+\mu'_2}{2},\frac{\mu'_1+\mu_2}{2})$ denote the joint measures whose marginals are  $\frac{\mu_1+\mu'_2}{2}$ and $\frac{\mu'_1+\mu_2}{2}$. Then for any joint measure $\mu\in \Gamma(\frac{\mu_1+\mu'_2}{2},\frac{\mu'_1+\mu_2}{2})$ we have
\[
\begin{split}
\int_{\Rn}u(y)\,d\big(\frac{\mu_1+\mu'_2}{2}\big)(y)\,-&\int_{\Rn}u(y)\,d\big(\frac{\mu'_1+\mu_2}{2}\big)(y)\\
=&\int_{\R^{2n}}\big(u(y)-u(y')\big)\,d\mu(y,y').
\end{split}
\]
Denoting
\[
\begin{split}
G(x,z):=u(x)-u(z),
\end{split}
\]
it holds that 
\begin{equation*}
\label{eq:dpp-2n}
\begin{split}
G(x,z)=&\sup_{\mu_1\in\mathcal{A}_1(x),\mu_2\in\mathcal{A}_2(z)}\bigg[\inf_{\mu'_2\in\mathcal{A}_2(x),\mu'_1\in\mathcal{A}_1(z)}\\&\bigg(\int_{\R^{2n}}G(y,y')\,d\mu_{\mu_1,\mu_1',\mu_2,\mu_2'}(y,y')\,\bigg)\bigg].
\end{split}
\end{equation*}
Different couplings generate different dynamic programming equations as well as  games in $\R^{2n}$,  and 
we may take the advantage of this flexibility in the proofs.

For example, with suitable couplings  a tug-of-war with noise described in  \cite{manfredipr12}  corresponds  to
\[
\begin{split}
G(x,z)=&\beta \int_{B(0,\eps)} G(x+h,y+P_{x,z}(h)) \ud h\\
&+ \frac\alpha2 \sup_{B(x,\eps)\times B(z,\eps)} G+\frac\alpha2 \inf_{B(x,\eps)\times B(z,\eps)} G,
\end{split}
\]
in $\R^{2n}$, where $P_{x,z}(h)$ is any isometry. 


Below we will prove  regularity results for several nonlinear problems of the form (\ref{DPP}). In all the examples of the paper, we require  that $\mu_i\in \mathcal A_i(x)$ implies $\spt \mu_i\in \ol B(x,\eps)$, where $\eps>0$ can be though as a step size in the underlying game. Throughout the paper, we denote by $B(x,\eps)$ an open ball centered at $x$ and of the radius $\eps>0$. When the  center point plays no role, we may drop it and denote $B_{\eps}$. 

\begin{theorem} 
\label{thm:main}
Let  $(x,z)\in B_R\times B_R$ and $B_{2R}\subset \Omega\,$. There exists $\delta=\delta(n)\in(0,1)$ such that if $u$ satisfies (\ref{DPP}) in the cases specified below, then
\begin{equation}
|u(x)-u(z)|\leq C\,\frac{|x-z|^{\delta}}{R^\delta}\,+\,C'(n)\frac{\eps^{\delta}}{R^\delta},
\end{equation}
where $C=C(n)(\sup_{B_{2R}\times B_{2R}}u(y)-u(y'))$. 
\end{theorem}

As shown in \cite{luirops13} this estimate can also be utilized in the proof of Harnack's inequality for the corresponding stochastic game. 

Without a loss of generality, we may assume that  $z=-x$, $R=1$ and $\sup_{B_2\times B_2}(u(y)-u(y'))\le 1$, by suitable translation, scaling and multiplication. Observe that  this does not require translation invariance.   
After the simplifications we see that it suffices to show that
\begin{equation}\label{modified}
|u(x)-u(-x)|\leq C |x|^{\delta}\,\,+\,C'(n)\eps^{\delta}\,,
\end{equation}
if $x\in B_1$ and $B_{2}\subset \Omega\,$. 


\subsection{Steps of the method}


As indicated in the introduction, our method  arises from game theory even if for expository reasons we have eventually avoided stochastic arguments. 
The first step of our regularity method is to  choose a comparison function $f:\R^{2n}\to \R$  that has the desired regularity properties. The key term in the comparison function we use to establish Hölder continuity in the examples below 
is $C \abs{x-z}^\delta$, $x,z\in \Rn$, $\delta\in(0,1)$.
 
The second step is to make a counter assumption: We wish to show that $u(x)-u(z)-f(x,z)$ is small enough in $B_1\times B_1\setminus T$, and in particular smaller than in $(B_2\times B_2)\setminus (B_1\times B_1\setminus T)$, and thrive for a contradiction by assuming that it is not.
    
 As a third step, we write down a multidimensional dynamic programming equation for the comparison function $f$ in $\R^{2n}$ by using a counter proposition.   
 
 As a final step, we derive a contradicting estimate for the dynamic programming equation of the previous step. Intuition based on game strategies gives guidelines how to obtain such an estimate by utilizing Taylor's expansion of $f(x,z)$ in suitable coordinates.

\subsection{Comparison function}
One of the key ideas of the proof is to use a suitably chosen \textit{comparison function}, defined in $\R^{2n}$.  We use throughout the work the function
\[
\begin{split}
f_1(x,z)=C|x-z|^{\delta}+|x+z|^2\,,
\end{split}
\]
where $C>1$ and
$\delta\in (0,1)$ are fixed later. The first term will give us the desired regularity estimate, and the second term makes sure that the estimate holds at  $(B_2\times B_2)\setminus (B_1\times B_1)$, see (\ref{vastaoletus}).

It is well known that functions satisfying the dynamic programming equation can be discontinuous at the small scale. Therefore, to prove the claim in the case $x$ and $z$ are close to each other, we have to use an auxiliary \textit{annular step function}. To this end, let $N=N(n)$ be an integer.
Then, set for $i=0,\ldots,N$ 
\[
\begin{split}
A_i=\{(x,z)\in\R^{2n}\,:\,(i-1)\frac{\eps}{10}<|x-z|\leq i\frac{\eps}{10}\}
\end{split}
\]
and define 
$f_2:\R^{2n}\to [0,\infty)$ by   
\begin{align}\label{stepfunction}
&f_2(x,z)=
\begin{cases}
 C^{2(N-i)}\eps^{\delta} & \text{ if } (x,z)\in A_i\,\text{, and }\\
 0 & \text{ if }|x-z|>N\frac{\eps}{10}.
\end{cases}
\end{align}
Observe that $f_2$ reaches its maximum $C^{2N}\eps^{\delta}$ on $A_0=\{(x,z)\,:\,x=z\}=T\,$. To see why we have chosen this $f_2$, one should consult   for example the calculations after (\ref{perusperus}).
 
Our final comparison function $f$ is composed of $f_1$ and $f_2$ as
\begin{equation}
f(x,z)=f_1(x,z)-f_2(x,z).
\end{equation}

We are going to produce estimates in terms of Taylor's expansion. It will be convenient to write the expansion in the form
\begin{equation}
\label{eq:taylor}
\begin{split}
f_1&(x+h_x,z+h_z)\\
&=f_1(x,z)+C\delta|x-z|^{\delta-1}(h_x-h_z)_V+2(x+z)\cdot (h_x+h_z)\\
&\,\,+\frac{C}{2}\delta|x-z|^{\delta-2}\big((\delta-1)(h_x-h_z)^2_V+(h_x-h_z)^2_{V^{\perp}}\big)\\
&\,\,+|h_x+h_z|^2\,+\mathcal{E}_{x,z}(h_x,h_z),
\end{split}
\end{equation}
where $V$ is the space spanned by $x-z$, $(h_x-h_z)_V$  refers to the scalar projection onto $V$ i.e.\ $(h_x-h_z)\cdot (x-z)/\abs{x-z}$, and $(h_x-h_z)_{V^{\perp}}$ onto the orthogonal complement. 

 The error term satisfies
\begin{equation}
|\mathcal{E}_{x,z}(h_x,h_z)|\leq C|(h_x,h_z)|^3(|x-z|-2\eps)^{\delta-3}\,, 
\end{equation}
if $|x-z|> 2\eps\,$. Especially, if we choose 
\begin{equation}\label{Noletus}
N\ge \frac{100C}{\delta}\,,
\end{equation}
then in the case $\abs{x-z}> N\frac{\eps}{10} $  and $|h_x|,|h_z|\leq \eps$, we can estimate
\begin{equation}
\label{error2}
\begin{split}
|\mathcal{E}_{x,z}(h_x,h_z)|&\leq C(2\eps)^3(\frac{|x-z|}{2})^{\delta-3}\,
\\
&\leq 
64C\eps^2|x-z|^{\delta-2}\frac{\eps}{|x-z|}\\
&\leq 64\eps^2|x-z|^{\delta-2}\frac{\delta}{10}\\
&\leq 10\eps^2|x-z|^{\delta-2}.
\end{split}
\end{equation}


\section{Tug-of-war}
\label{sec:p=infty}

We start with an example: we show how to obtain the asymptotic Hölder continuity for the value functions of the tug-of-war game by using the method of this paper.  There would be simpler methods available for the tug-of-war or the random walk, but the work done here will pay-off later, as we will see. Indeed, we  utilize these estimates later in  more difficult cases of Sections \ref{sec:p(x)} and \ref{sec:1<p}.  

The tug-of-war is a two-player zero-sum stochastic game played in $\Om$ with the following rules. Fix  the step size $\eps>0$ and an initial position $x_0\in
\Om$.  The players toss a fair coin and the winner of the
toss  moves the game position to any $x_1\in B_\eps
(x_0)$. The players continue playing
until the game position leaves the domain $\Omega$. At the end of the game
Player II pays Player I the amount determined  by a pay-off function $F$  defined outside $\Omega$.
Naturally, Player I tries to maximize the outcome while  Player II tries to minimize it, and thus the value of the game (for Player I) is defined as 
\[
u(x_0)=\sup_{S_{\I}}\inf_{S_{\II}}\,\mathbb{E}_{S_{\I},S_{\II}}^{x_0}[F(x_\tau)].
\]
Above $S_\I$ and $S_\II$ denote the strategies of Player I and Player II, respectively. If the players use strategies for which the game does not end almost surely and thus the expectation above is not well defined, we always set it to be $-\infty$.  

  The value function of the tug-of-war game satisfies the dynamic programming principle 
\begin{equation}
\label{eq:dpp-tgw}
\begin{split}
u(x)=\half\Big\{\sup_{B(x,\eps)}u+\inf_{B(x,\eps)} u\Big\}
\end{split}
\end{equation}
for $x\in \Om$, see \cite{peresssw09}  and \cite{lius15}. Intuitively, we obtain the value at the point $x$ by summing up the two possible outcomes of the coin toss with corresponding probabilities. The above equation is of the form (\ref{DPP}) and we can follow the steps introduced above.
 
\subsection{Multidimensional dynamic programming}
Because $u(x)-u(z)\le 1$ in $B_2\times B_2$ and $u(x)-u(z)=0$ on $T$, we see that  
\begin{equation}
\label{eq:boundary-bound}
\begin{split}
u(x)-u(z)-f(x,z)&= u(x)-u(z)-f_1(x,z)+f_2(x,z)\\
&\le \max f_2=C^{2N}\eps^{\delta}\,,
\end{split}
\end{equation}
if $(x,z)$ lies in $T$ or in  $(B_2\times B_2)\setminus (B_1\times B_1)$. We also used the fact that $f_1\ge 1$ in $(B_2\times B_2)\setminus (B_1\times B_1)$.

Next we will show that this inequality also holds in $B_1\times B_1\setminus T$. This would actually yield (\ref{modified}) with $C'(n)=C^{2N}\,$ and, eventually, the whole theorem. For the proof by contradiction, assume that
\begin{equation}\label{vastaoletus}
M:=\sup_{(x',z')\in B_1\times B_1\setminus T}(u(x')-u(z')-f(x',z'))\,>C^{2N}\eps^{\delta}\,. 
\end{equation}
Let then $\eta>0$. We choose $(x,z)\in B_1\times B_1\setminus T$ such that
\begin{align*}
u(x)-u(z)-f(x,z)\geq M-\eta\,.
\end{align*}

Then we are ready to write our first estimate:  by using (\ref{eq:dpp-tgw}) and the choices above, we obtain
\begin{equation}
\label{eq:tgw-main-est}
\begin{split}
M&\leq u(x)-u(z)-f(x,z)+\eta\\ 
&=\frac{1}{2}\sup_{B(x,\eps)}u
+\frac{1}{2}\inf_{B(x,\eps)}u
-\bigg(\frac{1}{2}\sup_{B(z,\eps)}u
+\frac{1}{2}\inf_{B(z,\eps)}u\bigg)
-f(x,z)+\eta\\
&=
\frac{1}{2}\bigg(\sup_{B(x,\eps)}u-\inf_{B(z,\eps)}u+\inf_{B(x,\eps)}u-\sup_{B(z,\eps)}u\,\bigg)
-f(x,z)+\eta\\
&=\,I_1-f(x,z)+\eta\,.
\end{split}
\end{equation}

We estimate $I_1$, and first observe that 
\begin{align}
\label{eq:sup-inf}
\sup_{B(x,\eps)}u-\inf_{B(z,\eps)}u\,\leq M+\sup_{B(x,\eps)\times B(z,\eps)}f+\eta\,.
\end{align}
This estimate follows easily by the above definition for $M$ as follows:
Suppose that $x_0\in B(x,\eps),z_0\in B(z,\eps)$ such that $\sup_{B(x,\eps)}u\le u(x_0)+\eta/2$ and $\inf_{B(z,\eps)} u\ge u(z_0)-\eta/2$. Then
\begin{align*}
\sup_{B(x,\eps)}u-\inf_{B(z,\eps)}u\,\le &u(x_0)-u(z_0)+\eta\\
=&u(x_0)-u(z_0)-f(x_0,z_0)+f(x_0,z_0)+\eta\\
\leq &M+f(x_0,z_0)+\eta\\
\leq& M+\sup_{B(x,\eps)\times B(z,\eps)}f+\eta\,.
\end{align*}
Observe that here we used the counter assumption. Indeed, the points $x_0$ and $z_0$ may lie outside $B_1\times B_1\setminus T$, but from (\ref{vastaoletus}) and the preceding discussion, we deduce that $M$ gives us the upper bound.

 Furthermore, analogous reasoning using (\ref{vastaoletus}) gives us that
\begin{align}
\label{eq:inf-sup}
\inf_{B(x,\eps)}u-\sup_{B(z,\eps)}u\,\leq M+\inf_{B(x,\eps)\times B(z,\eps)}f+\eta\,.
\end{align}
\sloppy
To be more precise,
choose $x_0\in B(x,\eps)$, $z_0\in B(z,\eps)$ such that $\inf_{B(x,\eps)\times B(z,\eps)}f \ge  f(x_0,z_0)-\eta$. Then
\begin{align*}
\inf_{B(x,\eps)}u-\sup_{B(z,\eps)}u&\leq  u(x_0)-u(z_0)\\
&=\, u(x_0)-u(z_0)-f(x_0,z_0)+f(x_0,z_0)\\
&\leq  \,M+f(x_0,z_0)\,\\
&\le  M+\inf_{B(x,\eps)\times B(z,\eps)}f+\eta\,
\end{align*}
and we obtain the above estimate.

Combining the estimates, we get
\[
\begin{split}
I_1 \leq M+\half \bigg(\sup_{B(x,\eps)\times B(z,\eps)}f \,+\,\inf_{B(x,\eps)\times B(z,\eps)}f\,\bigg)+\eta.
\end{split}
\]
From this and (\ref{eq:tgw-main-est}), we conclude that for the desired contradiction it suffices to show that 
\begin{equation}
\label{eq:f-p=infty}
\begin{split}
f(x,z)&>\frac{1}{2}\bigg(\sup_{B(x,\eps)\times B(z,\eps)}f \,+\,\inf_{B(x,\eps)\times B(z,\eps)}f\,\bigg).
\end{split}
\end{equation}

%

\subsection{Estimates}

Next we show the validity of the inequality (\ref{eq:f-p=infty}) for $f=f_1-f_2$, and for this we assume 
\begin{equation}\label{Coletus}
C=\frac{10^{10}}{\delta^{2}\omega}\,,
\end{equation}
for the constant $C$ in $f_1(x,z)=C\abs{x-z}^\delta+\abs{x+z}^2$. The constant $\om\in (0,1)$ will be later defined in this section. The choice is designed to cover also Sections \ref{sec:p=2} and \ref{sec:p(x)} with additional $n$-dependence $\om=\om(n)$. 
  If $|x-z|$ is large enough,    the validity of the inequality (\ref{eq:f-p=infty}) turns out to follow in a straightforward manner by using the Taylor expansion for $f_1$. The intuition coming from the underlying stochastic game will be useful in deriving the estimates. 

In the case  $|x-z|\approx \eps$, we will have to take the small scale jumps into account, and use the properties of the annular step function $f_2$ defined in (\ref{stepfunction}).

\subsubsection*{\bf Proof of inequality \eqref{eq:f-p=infty}, case $|x-z|> N\frac{\eps}{10}$}
Observe that the assumption $|x-z|> N\frac{\eps}{10}$ implies that $f_2(x,z)=0$, and therefore, it suffices to show (\ref{eq:f-p=infty}) for $f_1$.

Let $\eta >0$ and choose $h_x,h_z\in B(0,\eps)$ such that 
\[
\begin{split}
\sup_{B(x,\eps)\times B(z,\eps)}f_1\leq f_1(x+h_x,z+h_z)\,+\eta. 
\end{split}
\]
We first let $\theta=1/10$ and assume that
\begin{equation}\label{vastaan}
(h_x-h_z)_V^2\geq (4-\theta)\eps^2\,\,.
\end{equation}
Since $|h_x-h_z|< 2\eps$, this also means that 
\[
\begin{split}
(h_x-h_z)^2_{V^{\perp}}< \theta\eps^2\,.
\end{split}
\]
This can be though as fixing a strategy in the game in $\R^{2n}$. 

We get
by using the Taylor formula (\ref{eq:taylor}) and the estimate (\ref{error2}) for the error term that
\begin{align*}
&\sup_{B(x,\eps)\times B(z,\eps)}f_1 \,+\,\inf_{B(x,\eps)\times B(z,\eps)}f_1\,-2f_1(x,z)\\
&\le f_1(x+h_x,z+h_z)+f_1(x-h_x,z-h_z)-2f_1(x,z)+\eta\\
&=\, \frac{C}{2}\delta|x-z|^{\delta-2}\big(\,2(\delta-1)(h_x-h_z)_{V}^2+2(h_x-h_z)^2_{V^{\perp}}   \,\big)\\
&\,\,\,\,\,\,+2|h_x+h_z|^2+\mathcal{E}_{x,z}(h_x,h_z)+\mathcal{E}_{x,z}(-h_x,-h_z)+\eta\,\\
&\leq \frac{C}{2}\delta|x-z|^{\delta-2}\big(\,2(\delta-1)(4-\theta)\eps^2+2\theta\eps^2\,\big)+2(2\eps)^2\,+20\eps^2|x-z|^{\delta-2}+\eta\\
&\leq  |x-z|^{\delta-2}(20-C\delta)\eps^2+8\eps^2+\eta.
\end{align*}
The final inequality above follows e.g.\ by the choices $\theta=1/10$ and  $\delta\le1/10$. By using the assumption (\ref{Coletus}) for $C$, one observes that
\begin{equation}
\label{eq:vastaan-tulos}
|x-z|^{\delta-2}(20-C\delta)< -|x-z|^{\delta-2}10^8< -10^7\,.
\end{equation}
Combining this with above estimates, we conclude that
\[
\begin{split}
\sup_{B(x,\eps)\times B(z,\eps)}f_1+\inf_{B(x,\eps)\times B(z,\eps)}f_1-2f_1(x,z)<-10^6\eps^2\,.
\end{split}
\]

In turn, if (\ref{vastaan}) above does not hold, implying that
\begin{equation}\label{eivastaan}
(h_x-h_z)_V \leq (2-\frac{\theta}{4})\eps\,\,,  
\end{equation}
the desired estimate can be readily derived from the first order terms in Taylor estimate. More precisely, if $|h_x|,|h_z|\leq\eps$, the share of the second order term and the error terms in Taylor estimate (\ref{eq:taylor}) can be roughly estimated by
\begin{equation}
\label{eq:error-and-second-rough}
\begin{split}
&\frac{C}{2}\delta|x-z|^{\delta-2}(2\eps)^2+ (2\eps)^2+10\eps^2|x-z|^{\delta-2}\,\\
\leq & \,3C\delta\eps^2|x-z|^{\delta-2}\,\leq 3C\delta\eps|x-z|^{\delta-1}\frac{\eps}{|x-z|}\\
\leq & 3C\delta\eps|x-z|^{\delta-1}\frac{\eps}{N\frac{\eps}{10}}\leq \frac{30C}{N}\delta|x-z|^{\delta-1}\eps\,< \delta^2|x-z|^{\delta-1}\eps\,.
\end{split}
\end{equation}
Above we used the assumption (\ref{Noletus}) for $N\,$ and at the last step the fact that $\delta \in (0,1)$.
This in connection with the Taylor formula and  (\ref{eivastaan}) gives
\begin{align*}
&\sup_{B(x,\eps)\times B(z,\eps)}f_1+\inf_{B(x,\eps)\times B(z,\eps)}f_1-2f_1(x,z)\\
\leq & f_1(x+h_x,z+h_z)+f_1(x-\eps\frac{x-z}{|x-z|},z+\eps\frac{x-z}{|x-z|})-2f_1(x,z)\,+\eta\\
\leq &C\delta|x-z|^{\delta-1}((h_x-h_z)_V-2\eps)+4\abs{x+z}2\eps+\delta^2|x-z|^{\delta-1}\eps\,+\eta \,\\
\leq & C\delta|x-z|^{\delta-1}(-\frac{\theta}{4}\eps)+16\eps+\delta^2|x-z|^{\delta-1}\eps\,+\eta\\
\leq & (\delta-\frac{\theta}{4} C)\delta|x-z|^{\delta-1}\eps+16\eps +\eta.
\end{align*}
We point out that even if $\eps\frac{x-z}{|x-z|}\notin B(0,\eps)$, the estimate for $\inf$ above is valid due to continuity of $f_1$. Again, by applying the assumption (\ref{Coletus}) for $C$ and choosing $\om\le 1/10$, it is easy to check the desired negativity for this quantity. This completes the proof of the case $|x-z|> N\frac{\eps}{10}$. 

\subsubsection*{\bf Proof of inequality \eqref{eq:f-p=infty}, case $|x-z|\le N\frac{\eps}{10}$}  
In this case we can use the following elementary estimate for $f_1$:
If $x,z\in B_1$, and $h_x,h_z\in B(0,\eps),\,\eps<1$, then 
\begin{align}\label{perusperus}
|f_1(x+h_x,z+h_z)-f_1(x,z)|\leq 2C\eps^{\delta}+16\eps\,\leq 3C\eps^{\delta}\,
\end{align}
which immediately follows by using  the concavity and convexity of the terms in $f_1$.
We further obtain
\begin{equation}
\label{eq:tilde f}
\begin{split}
\sup_{B(x,\eps)\times B(z,\eps)}f_1-f_1(x,z)\leq  3C\eps^{\delta}\,.
\end{split}
\end{equation}
Moreover, 
\[
\begin{split}
\sup_{B(x,\eps)\times B(z,\eps)} (f_1-f_2)\le \sup_{B(x,\eps)\times B(z,\eps)} f_1-0.
\end{split}
\]
Since $(i-1)\frac{\eps}{10}< |x-z|\le i\frac{\eps}{10}$ for some $i=0,\ldots,N$,   we have 
\[
\begin{split}
\inf_{B(x,\eps)\times B(z,\eps)}  (f_1-f_2)
&\le \sup_{B(x,\eps)\times B(z,\eps)} f_1-\sup_{B(x,\eps)\times B(z,\eps)} f_2 \\
&\le \,\,\sup_{B(x,\eps)\times B(z,\eps)} f_1- C^{2(N-i+1)}\eps^{\delta}\\
&\le\sup_{B(x,\eps)\times B(z,\eps)} f_1-C^{2}C^{2(N-i)}\eps^{\delta}\\ 
&= \,\sup_{B(x,\eps)\times B(z,\eps)} f_1-(C^{2}-2)C^{2(N-i)}\eps^{\delta}-2C^{2(N-i)}\eps^{\delta}\,\\
&=\sup_{B(x,\eps)\times B(z,\eps)} f_1-\, (C^{2}-2)C^{2(N-i)}\eps^{\delta}-2f_2(x,z)\,\\
&\le \sup_{B(x,\eps)\times B(z,\eps)} f_1-10C\eps^{\delta}-2f_2(x,z)\,.
\end{split}
\]
Summing up the previous two estimates and recalling (\ref{eq:tilde f}), we end up with 
\[
\begin{split}
\sup_{B(x,\eps)\times B(z,\eps)}f+&\inf_{B(x,\eps)\times B(z,\eps)} f\\
&\le \sup_{B(x,\eps)\times B(z,\eps)} f_1+\sup_{B(x,\eps)\times B(z,\eps)} f_1 -10C\eps^\delta-2f_2(x,z)\\
&\le 2 f_1(x,z)+6C\eps^\delta  -10C\eps^\delta-2f_2(x,z)\\
&< 2 f(x,z). 
\end{split}
\]
This proves the estimate \eqref{eq:f-p=infty}, and thus (\ref{modified}) and finally Theorem \ref{thm:main} in the case of the tug-of-war.

\begin{remark}
\label{rem:inf-harm}
By passing to a limit this implies Hölder continuity for the infinity harmonic functions, cf.\ \cite{peresssw09} or the proof of Theorem  4.9 in \cite{manfredipr12}, i.e.\ for the viscosity solutions to 
\[
\begin{split}
\Delta_\infty u=\sum_{i,j=1}^n u_{ij} u_i u_j=0,
\end{split}
\]
where $u_i$ and $u_{ij}$ denote the elements of the gradient and the Hessian, respectively.  As already pointed out, this is not our main motivation but we use the estimates of this section as a tool in Section \ref{sec:p(x)} without repeating them there.
\end{remark}

\section{Random walk}
\label{sec:p=2}

We continue with another well known example: we show the asymptotic Hölder continuity for the random walk. Again some estimates in more difficult cases in Sections \ref{sec:p(x)} and \ref{sec:1<p} are essentially the same, and we do not need to repeat the estimates there. Therefore this expository choice does not add much to the length of the paper. 

We consider a version of the random walk where at $x\in \Om$ the next point is chosen according to the uniform probability distribution on $B(x,\eps)$, and this is repeated until the process  exits $\Om$. At this time
 the amount given by a pay-off function $F$  defined outside $\Omega$ is collected.
The function $u$ is the expected pay-off  of this process
\[
u(x_0)=\mathbb{E}^{x_0}[F(x_\tau)].
\]
This expectation can also be written as an average over the neighboring expected pay-offs as\begin{equation}
\label{eq:dpp-random-walk}
\begin{split}
u(x)=\kint_{B(x,\eps)} u\ud y:=\frac1{\abs{B(x,\eps)}} \int_{B(x,\eps)} u \ud y.
\end{split}
\end{equation}
This again is of the form (\ref{DPP}) and we can employ the same method as above.

\subsection{Multidimensional dynamic programming}
Similarly as in Section \ref{sec:p=infty}, the equation
 (\ref{eq:dpp-random-walk}) gives
\begin{align*}
M&\leq u(x)-u(z)-f(x,z)+\eta\\ 
&=\bigg(\intav_{B(x,\eps)}u(y)\,dy\,-\intav_{B(z,\eps)}u(y)\,dy\,\bigg)
-f(x,z)+\eta\\
&=\,I_2-f(x,z)+\eta\,.
\end{align*}

Define $P_{x,z}(h)$ as a mirror point of $h$ with respect to $V^\perp=\operatorname{span}(x-z)^\perp$.
The term $I_2$ (we later utilize  this estimate in Section \ref{sec:p(x)}, which explains the notation) can be written as
\begin{align*}
I_2&=\frac{1}{|B_{\eps}|}\bigg(\int_{B(x,\eps)\setminus B(z,\eps)}u(y)\,dy-\int_{B(z,\eps)\setminus B(x,\eps)}u(y)\,dy\\
&\,\,\,\,\,\,+\int_{B(x,\eps)\cap B(z,\eps)}u(y)-u(y)\,dy\,\bigg)\\
&=\frac{1}{|B_{\eps}|}\bigg(\int_{B(0,\eps)\setminus B(z-x,\eps)}u(x+h)-u(z+P_{x,z}(h))\,dh\\
&\,\,\,\,\,\,+\int_{B(x,\eps)\cap B(z,\eps)}0\,dy\,\bigg).
\end{align*}
Then by adding and subtracting we obtain
\begin{align*}
&I_2
=\frac{1}{|B_{\eps}|}\bigg(\int_{B(0,\eps)\setminus B(z-x,\eps)}\Big(u(x+h)-u(z+P_{x,z}(h))\\
&\hspace{15em}-f(x+h,z+P_{x,z}(h))\Big)\,dh\\
&\,\,\,\,\,\,-\int_{B(x,\eps)\cap B(z,\eps)}f(y,y)\,dy\,\\
&\,\,\,\,\,\,+\int_{B(0,\eps)\setminus B(z-x,\eps)}f(x+h,z+P_{x,z}(h))\,dh\,+\int_{B(x,\eps)\cap B(z,\eps)}f(y,y)\,dy\,\bigg)\\
&\leq M+\frac{1}{|B_{\eps}|}\bigg(\,\int_{B(0,\eps)\setminus B(z-x,\eps)}f(x+h,z+P_{x,z}(h))\,dh\,
+\int_{B(x,\eps)\cap B(z,\eps)}f(y,y)\,dy\,\bigg)\,,
\end{align*}
where we again used the definition of $M$, as well as the counter assumption and the definition of $f_2$ to estimate $-f(y,y)=f_2(y,y)\le M$.

From the previous estimate we conclude the desired contradiction if at every $(x,z)\in B_1\times B_1\setminus T$ it holds that
\begin{equation}
\label{eq:f-p=2}
\begin{split}
f(x,z)>\frac{1}{|B_{\eps}|}\bigg(\,\int_{B(0,\eps)\setminus B(z-x,\eps)}&f(x+h,z+P_{x,z}(h))\,dh\,\\
&\hspace{1 em}+\int_{B(x,\eps)\cap B(z,\eps)}f(y,y)\,dy\,\bigg)\,.
\end{split}
\end{equation}

\subsection{Estimates}
Similarly as in Section \ref{sec:p=infty}, we estimate the two cases separately.
\subsubsection*{\bf Inequality \eqref{eq:f-p=2}, case $|x-z|> N\frac{\eps}{10}$}
In this case the claim reduces to 
\begin{align*}
f_1(x,z)&>\frac{1}{|B_{\eps}|}\int_{B(0,\eps)}f_1(x+h,z+P_{x,z}(h))\,dh\,.
\end{align*}
We utilize the Taylor series (\ref{eq:taylor}). When integrating the series, the first order terms vanish by symmetry, and by the definition of the mirror mapping $P_{x,z}$ also $(h-P_{x,z}(h))^2_{V^\perp}=0$. This and the error estimate (\ref{error2}), gives 
\begin{align*}
&\frac{1}{|B_{\eps}|}\int_{B(0,\eps)}f_1(x+h,z+P_{x,z}(h))\,dh-f_1(x,z)\\
&=\,\frac{2C\delta(\delta-1)|x-z|^{\delta-2}}{|B_{\eps}|}\int_{B(0,\eps)}
h_V^2\,dh\,+\frac{1}{|B_{\eps}|}\int_{B(0,\eps)}\mathcal{E}_{x,z}(h,P_{x,z}(h))\,dh\,\\
&\leq \,2C\delta(\delta-1)|x-z|^{\delta-2} \frac{\eps^2}{n+2}
+10\eps^2|x-z|^{\delta-2}\,\\
&\leq\,\eps^2|x-z|^{\delta-2}\big(10-\frac{C\delta}{4(n+2)}\big)\,,
\end{align*}
where by a direct calculation $\kint_{B(0,\eps)} h_V^2 \ud y=\eps^2/(n+2)$ and we assumed $\delta\le1/10$. Observe that since  $P_{x,z}(h)$ gives the mirror point of $h$, the terms of the form  $(\cdot)_{V^\perp}$ vanished. Moreover, by symmetry the first order terms vanished as well.
Now, the assumption (\ref{Coletus}) on $C$ with $\om=\om(n)\le 1/(n+2)$  guarantees the negativity of $\big(10-\frac{C\delta}{4(n+2l)}\big)\,$. This completes the proof of the present case.
\subsubsection*{\bf Inequality \eqref{eq:f-p=2}, case $|x-z|\le N\frac{\eps}{10}$}
For this, suppose first that $|x-z|\geq \frac{7}{4}\eps$. In this case, observe that 
\begin{align*}
B\bigg(\frac{\eps(z-x)}{2|z-x|},\frac{\eps}{4}\bigg)\subset B(0,\eps)\setminus B(z-x,\eps)\,,
\end{align*}
and  
\begin{align*}
|x+h-(z+P_{x,z}(h))|\leq |x-z|-\frac{\eps}{10}
\end{align*}
for every $h\in B(\frac{\eps(z-x)}{2|z-x|},\frac{\eps}{4}).$ 
This in turn implies, by the definition of the annular function $f_2$, that 
\[
\begin{split}
f_2(x+h,z+P_{x,z}(h))\geq C^2f_2(x,z)\,\,\text{ if }h\in B\big(\frac{\eps(z-x)}{2|z-x|},\frac{\eps}{4}\big)\,.
\end{split}
\]
This implies
\begin{align*}
&\frac{1}{|B_{\eps}|}\bigg(\,\int_{B(0,\eps)\setminus B(z-x,\eps)}f_2(x+h,z+P_{x,z}(h))\,dh\,+\int_{B(x,\eps)\cap B(z,\eps)}f_2(y,y)\,dy\,\bigg)\,\\
\,&>\frac{1}{|B_{\eps}|}\,\int_{B(\frac{\eps(z-x)}{2|z-x|},\frac{\eps}{4})}f_2(x+h,z+P_{x,z}(h))\,dh\\
\,&>\frac{|B_{\frac{\eps}{4}}|}{|B_{\eps}|}C^2f_2(x,z)=\,\frac{C^2}{4^n}f_2(x,z)\,.
\end{align*}
Furthermore, as before we use the rough estimate (\ref{perusperus}) for  $f_1$, implying that the inequality (\ref{eq:f-p=2}) can fail for $f_1$ at most by $3C\eps^{\delta}\,$. Combining this with the above estimate for $f_2$, we conclude that
\begin{align*}
\frac{1}{|B_{\eps}|}&\bigg(\,\int_{B(0,\eps)\setminus B(z-x,\eps)}f(x+h,z+P_{x,z}(h))\,dh\,+\int_{B(x,\eps)\cap B(z,\eps)}f(y,y)\,dy\,\bigg)\,\\
&< f_1(x,z)+\,3C\eps^{\delta}-\frac{C^2}{4^n}f_2(x,z)\,\\
&<\,f_1(x,z)-f_2(x,z)+\,3C\eps^{\delta}-(\frac{C^2}{4^n}-1)f_2(x,z)<f(x,z)\,,
\end{align*}
where the final inequality follows simply by using again the size condition (\ref{Coletus}) for $C\,$ and choosing $\om\le4^{-n}$. This completes the proof of the case $|x-z|\geq \frac{7}{4}\eps$. 

For the final case $|x-z|<\frac{7}{4}\eps$, observe that
\begin{align*}
|B(x,\eps)\cap B(z,\eps)|>\frac{1}{4^n}|B_{\eps}|\,. 
\end{align*}
Moreover, because $x\neq z$ 
\[
\begin{split}
f_2(y,y)\geq C^2f_2(x,z)\,.
\end{split}
\]
Then we simply estimate that 
\begin{align}
\label{eq:weird-term-used}
\frac{1}{|B_{\eps}|}&\,\int_{B(0,\eps)\setminus B(z-x,\eps)}f_2(x+h,z+P_{x,z}(h))\,dh\,\notag\\
&\hspace{5 em}+\frac{1}{|B_{\eps}|}\int_{B(x,\eps)\cap B(z,\eps)}f_2(y,y)\,dy\, \\
&>\frac{1}{|B_{\eps}|}\int_{B(x,\eps)\cap B(z,\eps)}f_2(y,y)\,dy\,\notag
\\
&>\frac{C^2}{4^n}f_2(x,z)\,.\notag
\end{align} 
Then the desired inequality for $f=f_1-f_2$ follows exactly in the same way as in the case $|x-z|\geq\frac{7}{4}\eps\,$.
Again Theorem \ref{thm:main} follows.

As an application, by passing to a limit, this implies local Hölder continuity for harmonic functions similarly as explained for infinity harmonic functions in Remark \ref{rem:inf-harm}. By adjusting the comparison function, it would also be possible to obtain local Lipschitz continuity. However, for the consistency of the exposition we shall not pursue this analysis here.


\section{Tug-of-war with noise and space dependent probabilities}

\label{sec:p(x)}

Next we proceed to new results. We consider a tug-of-war with noise and space dependent probabilities:
This is again a two-player zero-sum stochastic game played in a domain $\Om$. The rules are as
follows. Fix  the step size $\eps>0$ and an initial position $x_0\in
\Om$. 
 The players start by tossing a biased space dependent coin with probabilities
$\alpha(x_0)$ and $\beta(x_0)$, $\alpha(x_0) + \beta(x_0) =1$. If the result is heads
(with probability $\alpha(x_0)$), then they play a tug-of-war as described in Section \ref{sec:p=infty}. 
On the other hand, if the result  is tails (with probability $\beta(x_0)$), the game
state moves to a random point
in the ball $B(x_0,\eps)$. The players continue playing
until the game position leaves the domain $\Omega$. 
Then Player II pays Player I the amount determined  by a pay-off function $F$  defined outside $\Omega$.

Similarly as in the tug-of-war, Player I tries to maximize the pay-off while  Player II tries to minimize it and thus \emph{value of the game} is defined as  
\[
u(x_0)=\sup_{S_{\I}}\inf_{S_{\II}}\,\mathbb{E}_{S_{\I},S_{\II}}^{x_0}[F(x_\tau)],
\]
where $S_\I$ and $S_\II$ again denote the strategies of Player I and Player II, respectively. The usual tug-of-war with noise with fixed probabilities in \cite{manfredipr12}, the tug-of-war in Section \ref{sec:p=infty}, and the random walk in Section \ref{sec:p=2} are special cases of this game. 

Deriving the dynamic programming principle  for this game is beyond the scope of the paper (for the fixed probabilities case, see \cite{luirops13}). Nonetheless, we directly take the dynamic  programming  equation
\begin{equation}
\label{eq:dpp-p(x)}
\begin{split}
u(x) = \frac{\alpha(x)}{2} \left\{ \sup_{ B(x,\eps)} u +
\inf_{ B(x,\eps)} u \right\} + \beta(x) \kint_{ B(x,\eps)}
u \ud y
\end{split}
\end{equation}
as our starting point. Heuristically we obtain the value at the point $x\in \Om$ by summing up the three possible outcomes with corresponding probabilities. This equation is again of such a form that we can apply our method. We also point out that we have no regularity assumption for $\alpha(x)$.

\subsection{Multidimensional dynamic programming}

Let $\eps>0$ and  $\alpha:\Omega\to[0,1]$ be a Borel measurable mapping.
Without loss of generality, we may assume $\max\{\alpha(x),\alpha(z)\}=\alpha(x)$ and 
$u(x)-u(z)\ge 0$. Indeed, if $u(x)-u(z)< 0$, then consider $-u$ instead. Along the same lines as before, the equation  (\ref{eq:dpp-p(x)}) gives  
\begin{align*}
M\leq& u(x)-u(z)-f(x,z)+\eta\\ 
=&\frac{\alpha(x)}{2}\sup_{B(x,\eps)}u
+\frac{\alpha(x)}{2}\inf_{B(x,\eps)}u
+\,\beta(x)\intav_{B(x,\eps)}u(y)\,dy\,\\
&-\frac{\alpha(z)}{2}\sup_{B(z,\eps)}u
-\frac{\alpha(z)}{2}\inf_{B(z,\eps)}u
-\,\beta(z)\intav_{B(z,\eps)}u(y)\,dy\,-f(x,z)+\eta\\
=&
\frac{\alpha(z)}{2}\bigg(\sup_{B(x,\eps)}u-\inf_{B(z,\eps)}u+\inf_{B(x,\eps)}u-\sup_{B(z,\eps)}u\,\bigg)
\\
&+\,\beta(x)\bigg(\intav_{B(x,\eps)}u(y)\,dy\,-\intav_{B(z,\eps)}u(y)\,dy\,\bigg)
\\
&+\frac{(\alpha(x)-\alpha(z))}{2}\bigg(\sup_{B(x,\eps)}u+\inf_{B(x,\eps)}u-2\intav_{B(z,\eps)}u(y)\,dy\,\bigg)-f(x,z)+\eta\\
=&\,I_1+I_2+I_3-f(x,z)+\eta\,.
\end{align*}

The terms $I_1$ and $I_2$ were already estimated in the previous sections, so we may focus our attention on $I_3$.
Choose a sequence $x_k$ such that $\lim_k u(x_k)=\sup_{B(x,\eps)}u
$. We have 
\begin{align*}
\sup_{B(x,\eps)}&u-\intav_{B(z,\eps)}u(y)\,dy\,=
\intav_{B(z,\eps)} \lim_k u(x_k)-u(y)\,dy\,\\
=&\intav_{B(z,\eps)}\lim_k (u(x_k)-u(y)-f(x_k,y)+f(x_k,y))\,dy\,\\
\leq & M+\sup_{x'\in B(x,\eps)} \intav_{B(z,\eps)}f(x',y)\,dy\,.
\end{align*}
By a reversed reasoning, one can also verify that
\begin{align*}
\inf_{B(x,\eps)}u-\intav_{B(z,\eps)}u(y)\,dy\,\leq
M+\intav_{B(z,\eps)}\inf_{\tilde{x}\in B(x,\eps)}f(\tilde{x},y)\,dy\,.
\end{align*}
To see this, write 
\begin{align*}
&\inf_{B(x,\eps)}u-\intav_{B(z,\eps)}u(y)\,dy\,=\intav_{B(z,\eps)}\inf_{B(x,\eps)}u-u(y)\,dy\,\\
&= \intav_{B(z,\eps)}\inf_{B(x,\eps)}u-u(y)-\inf_{\tilde{x}\in B(x,\eps)}f(\tilde{x},y)\,dy\,+
\intav_{B(z,\eps)}\inf_{\tilde{x}\in B(x,\eps)}f(\tilde{x},y)\,dy\,
\\
&\leq M+\,\intav_{B(z,\eps)}\inf_{\tilde{x}\in B(x,\eps)}f(\tilde{x},y)\,dy\,.
\end{align*}

By combining these estimates, we get that
\begin{align*}
I_3\leq \bigg(\frac{\alpha(x)-\alpha(z)}{2}\bigg)\bigg(2M+\sup_{x'\in B(x,\eps)}\intav_{B(z,\eps)}f(x',y)+\inf_{\tilde{x}\in B(x,\eps)}f(\tilde{x},y)\,dy\,\bigg)\,.
\end{align*}

Summing up and also recalling the estimates from the previous sections, to obtain the desired contradiction it suffices to show that $f$ satisfies the following three conditions at every $(x,z)\in B_1\times B_1\setminus T$:
\begin{align*}
\tag{I}f(x,z)&>\frac{1}{2}\bigg(\sup_{B(x,\eps)\times B(z,\eps)}f \,+\,\inf_{B(x,\eps)\times B(z,\eps)}f\,\bigg)\,,\\
\tag{II} f(x,z)&>\frac{1}{|B(x,\eps)|}\bigg(\,\int_{B(0,\eps)\setminus B(z-x,\eps)}f(x+h,z+P_{x,z}(h))\,dh\,\\
&\hspace{10 em}+\int_{B(x,\eps)\cap B(z,\eps)}f(y,y)\,dy\,\bigg)\,,\\
\tag{III} f(x,z)&>\frac{1}{2}\sup_{x'\in B(x,\eps)}\bigg[\intav_{B(z,\eps)}f(x',y)+\inf_{\tilde{x}\in B(x,\eps)}f(\tilde{x},y)\,dy\,\bigg].
\end{align*}

 The inequalities I and II were treated earlier in Sections \ref{sec:p=infty} and \ref{sec:p=2}, and we may concentrate on the inequality III.

\subsection{Estimates}
\subsubsection*{\bf Inequality III, case $|x-z|>  N\frac{\eps}{10}$} 
As in the previous case, it suffices to verify the desired inequality for $f_1$. 
Let $\eta >0$. Suppose that $h_x\in B(0,\eps)$ is such that
\begin{equation}
\begin{split}
\label{eq:sup-direction}
\sup_{x'\in B(x,\eps)}&\bigg[\intav_{B(z,\eps)}f_1(x',y)+\inf_{\tilde{x}\in B(x,\eps)}f_1(\tilde{x},y)\,dy\,\bigg]\\
\leq\,&\intav_{B(z,\eps)}f_1(x+h_x,y)+\inf_{\tilde{x}\in B(x,\eps)}f_1(\tilde{x},y)\,dy\,\,+\eta\,=: I_3^1.
\end{split}
\end{equation}
Next we estimate $\inf_{\tilde{x}\in B(x,\eps)}f_1(\tilde{x},y)\le f_1(x-h_x,y)$ and use the Taylor expansion for $f_1$ to obtain
\begin{align*}
I_3^1&\leq  \intav_{B(z,\eps)}f_1(x+h_x,y)\,dy+f_1(x-h_x,y)\,dy\,+\eta\,\\
&=\intav_{B(0,\eps)}f_1(x+h_x,z+h_z)+f_1(x-h_x,z-h_z)\,dh_z\,+\eta\,\\
&= \intav_{B(0,\eps)}2f_1(x,z)+C\delta|x-z|^{\delta-2}\big((\delta-1)(h_x-h_z)_{V}^2+(h_x-h_z)^2_{V^{\perp}}\big)\,dh_z\,\\
&\,\,\,\,\,\,+\intav_{B(0,\eps)}2|h_x+h_z|^2+\mathcal{E}_{x,z}(h_x,h_z)+\mathcal{E}_{x,z}(-h_x,-h_z)\,dh_z\,+\eta\\
&\leq 2f_1(x,z)+C\delta|x-z|^{\delta-2}\intav_{B(0,\eps)}(\delta-1)(h_x-h_z)_{V}^2+(h_x-h_z)^2_{V^{\perp}}\,dh_z\\
&\hspace{1 em}+2(2\eps)^2+2(10\eps^2|x-z|^{\delta-2})+\eta.\end{align*}
First, we consider the case, where $h_x$ lies close to $\eps\frac{x-z}{|x-z|}$. To be more precise, let
\begin{equation}\label{ekaII}
\bigg|h_x-\frac{x-z}{|x-z|}\eps\bigg|\leq \theta(n)\eps\,
\end{equation}
such that the estimate
\[
\begin{split}
\intav_{B(0,\eps)}(\delta-1)(h_x-h_z)_{V}^2+(h_x-h_z)^2_{V^{\perp}}\,dh_z\le -\frac{1}{n+2}\eps^2.
\end{split}
\]
holds for any $\delta\leq \frac{1}{10(n+2)}$ and $\theta=\frac{1}{10(n+2)}$ by a direct estimation.
Indeed, first observe that
\[
\begin{split}
\intav_{B(0,\eps)}&-\Big(\frac{x-z}{|x-z|}\eps-h_z\Big)_{V}^2+\Big(\frac{x-z}{|x-z|}\eps-h_z\Big)^2_{V^{\perp}}\,dh_z
=-\eps^2+\frac{n-2}{n+2}\eps^2
\end{split}
\]
and then take the small adjustments into account.
Combining the above estimates, we obtain in the case (\ref{ekaII}) that
\[
\begin{split}
I_3^1\le 2f_1(x,z)-C\delta|x-z|^{\delta-2}\frac{\eps^2}{n+2}+8\eps^2+20\eps^2|x-z|^{\delta-2}\,+\eta.
\end{split}
\]
Using the assumption (\ref{Coletus}) on the size of $C$ we get
\begin{align*}\label{tapaus2}
8\eps^2&+20\eps^2|x-z|^{\delta-2}-C\delta|x-z|^{\delta-2}\frac{\eps^2}{n+2}+\eta\\
<\,&100\eps^2|x-z|^{\delta-2}-C\delta|x-z|^{\delta-2}\frac{\eps^2}{n+2}+\eta\\
<\,&\eps^2|x-z|^{\delta-2}(100-\frac{C\delta }{n+2})\\
<\,&-\eps^2\,,
\end{align*}
where we choose $\om$ in the bound for  $C$ so that $\om \le\frac{1}{n+2}$, and $\eta$ small enough. This verifies the claim in the case (\ref{ekaII}). 

In the remaining case
\begin{equation}
\label{eq:larger}
\bigg|h_x-\frac{x-z}{|x-z|}\eps\bigg|> \theta(n)\eps=\frac{1}{10(n+2)}\eps,
\end{equation}
we utilize the first order terms in the Taylor expansion.
To this end, by (\ref{eq:error-and-second-rough}), the second order and error terms are bounded by $\delta^2|x-z|^{\delta-1}\eps\,$. Recalling \eqref{eq:sup-direction} and estimating $\inf_{\tilde{x}\in B(x,\eps)}f_1(\tilde{x},y)\le f_1(x-\eps\frac{x-z}{|x-z|},y)$, we have 
\begin{align*}
&I_3^1 \le\intav_{B(0,\eps)}f_1(x+h_x,z+h_z)+f_1(x-\eps\frac{x-z}{|x-z|},z-h_z)\,dh_z\,+\eta\,\\
&\leq 2f_1(x,z)+ \intav_{B(0,\eps)}C\delta|x-z|^{\delta-1}\big((h_x-h_z)_{V}+(-\eps\frac{x-z}{|x-z|}-(-h_z))_V\big)\,dh_z\,\\
&\,\,\,+\intav_{B(0,\eps)}2(x+z)\cdot \big((h_x+h_z)_V+(-\eps\frac{x-z}{|x-z|}+(-h_z))_V\big)\,dh_z\,\\
&\,\,\,+\delta^2|x-z|^{\delta-1}\eps\,+\eta.
\end{align*}
In the first integral above, we estimate $(-h_z)_V+(h_z)_V=0$ and $(-\eps\frac{x-z}{|x-z|})_V=-\eps$. The second integral on the right hand side above is bounded by $16\eps$. Thus by using (\ref{eq:larger}), we obtain
\begin{align*}
I_3^1&\le 2f_1(x,z)+C\delta|x-z|^{\delta-1}((h_x)_V-\eps)+16\eps+\delta^3|x-z|^{\delta-1}\eps\,+\eta\\
&\leq 2f_1(x,z)-C\delta|x-z|^{\delta-1}\frac{\theta}{2} \eps +\delta^2|x-z|^{\delta-1}\eps\,+\eta\\
&\leq 2f_1(x,z)-\delta |x-z|^{\delta-1}\eps\big(\frac{C\theta}{2} -\delta\,\big)+\eta\\
&\leq 2f_1(x,z)-\eps\,.
\end{align*}
Again the final inequality follows recalling the bound (\ref{Coletus}) for $C$ and choosing $\om\le \frac{1}{10(n+2)}$. 
This completes the proof of the present case.

\subsubsection*{\bf Inequality III, case $|x-z|\le  N\frac{\eps}{10}$}
In this case we again  use   the rough estimate (\ref{perusperus}) for $f_1$ in a similar spirit as before. 

Consider then the quantity 
\begin{align*}
\bigg(\sup_{x'\in B(x,\eps)}\bigg[\intav_{B(z,\eps)}-f_2(x',y)+\inf_{\tilde{x}\in B(x,\eps)}-f_2(\tilde{x},y)\,dy\,\bigg]\bigg)+2f_2(x,z)
\end{align*}
for $f_2$. We use the term  $\intav_{B(z,\eps)}\inf_{\tilde{x}\in B(x,\eps)}-f_2(\tilde{x},y)\,dy$, and choose $\tilde x$ such that the distance between this choice and $y$ is suitably smaller than $\abs{x-y}$, so that we get the desired estimate.

For that, let us define for every $y\in B(z,\eps)$, $P_x(y)\in B(x,\frac{\eps}{2})$ such that
\[
\begin{split}
P_x(y)=
\begin{cases}
x\,\,&\text{ if }|x-y|\geq\frac{\eps}{2}\\
y\,\,&\text{ if }|x-y|<\frac{\eps}{2}\,.
\end{cases}
\end{split}
\]
Then for any $x'\in  B(x,\eps)$ it holds that
\begin{align*}
&\intav_{B(z,\eps)}-f_2(x',y)+\inf_{\tilde{x}\in B(x,\eps)}(-f_2(\tilde{x},y))\,dy\\
&\le \intav_{B(z,\eps)}-f_2(x',y)-f_2(P_x(y),y)\,dy\,
\le  -\intav_{B(z,\eps)}f_2(P_x(y),y)\,dy\,.
\end{align*}

Then suppose that  $|x-z|\geq\frac{2\eps}{3}$, denote $z'=z+\frac{2}{3}\eps\frac{x-z}{|x-z|}$ and $y\in B(z',\frac{\eps}{3})\subset B(z,\eps)\,$. It follows that 
\[
\begin{split}
|P_x(y)-y|&\leq \abs{x-y} \le \abs{x-z'}+\frac{\eps}{3}\le \abs{x-z}-\frac{\eps}{3}.
\end{split}
\]
Thus by the definition of $f_2$ and the fact that always $x\neq z$, it follows that
\begin{align*}
f_2(P_x(y),y)\geq C^2f_2(x,z)
\end{align*}
in the same ball.
Using this and the fact that $f_2\ge 0$, we obtain
\begin{align*}
\intav_{B(z,\eps)}&f_2(P_x(y),y)\,dy\,\\
&\ge \frac1{\abs{B_\eps}} \bigg(\int_{B(z,\eps)\cap B(z',\frac{\eps}{3})} f_2(P_x(y),y) \ud y+ \int_{B(z,\eps)\setminus B(z',\frac{\eps}{3})} f_2(P_x(y),y) \ud y \bigg)\\
&\ge \frac{C^2}{3^n}f_2(x,z)>2f_2(x,z)+6C\eps^{\delta}, 
\end{align*}
 where the above inequalities follows from (\ref{Coletus}), the choice $\om^2\le 3^{-n}$  and the definition of $f_2$.

In turn, if $|x-z|<\frac{2\eps}{3}$, then we use the fact that $P_x(y)=y$ in $B(x,\frac{\eps}{3})\subset B(z,\eps)$. This and the definition of $f_2$ imply  
\begin{align*}
\intav_{B(z,\eps)}f_2(P_x(y),y)\,dy&\ge \frac1{\abs{B_\eps}}\int_{B(x,\frac{\eps}{3})}f_2(P_x(y),y)\,dy\\
&\ge  \frac{C^2}{3^n}f_2(x,z) > 2f_2(x,z)+6C\eps^{\delta}.
\end{align*}
At the last step we again used (\ref{Coletus}) and properties of $f_2\,$. This completes the proof of Theorem \ref{thm:main}. Observe that the choices $\delta=\frac{1}{10(n+2)}$ and $\omega=\min\{\frac{1}{10(n+2)},\frac{1}{4^n}\}$ are sufficient with respect to all the choices made above.


The above result can be applied to the game theoretic approach to the normalized $p(x)$-Laplacian
\[
\begin{split}
\Delta_{p(x)}^N u=:\Delta u+(p(x)-2)\Delta_{\infty}^N u=0,\quad\text{for}\quad x\in \Om \Subset \Rn,  
\end{split}
\]
where $\Delta_{\infty}^N u=\abs{\nabla u}^{-2} \sum_{i,j=1}^n u_{ij} u_i u_j$ with the notation of Remark \ref{rem:inf-harm}, and $2\le p(x)\le \infty$. Indeed, the regularity estimate derived above allows us to pass to a limit, and thus to obtain estimates for viscosity solutions (see \cite{juutinenlp10}, as well as \cite{julin15} for a different version of the $p(x)$-Laplacian) to the above equation using game theory by setting 
\[
\begin{split}
\alpha(x)=\frac{p(x)-2}{p(x)+n},\qquad \beta(x)=\frac{2+n}{p(x)+n}.
\end{split}
\]

It is also worth noting that the approach in \cite{peress08} and \cite{manfredipr12} does not directly extend to the game of this section, since these papers utilize the translation invariance.

%
%


\section{Tug-of-war with noise}
\label{sec:1<p}

The previous sections covered the game related to $p$-harmonic functions in the case $2\le p\le \infty$. In this section, we will include the range $1<p<2$ by considering a slight modification of the game in \cite{peress08} suggested in \cite{kawohlmp12}.  

Choose a starting point $x_0\in
\Om$. The players fix their possible moves $\nu_\I$ and $\nu_\II$ with $\abs{\nu_I},\abs{\nu_\II}\leq \eps$, and toss a fair coin. If Player I wins the toss, then she tosses a biased coin. If she gets heads (with probability $\alpha>0$), the token is placed at $x_1=x_0+\nu_\I$. If she gets tails (with probability $\beta$), then the token is placed at a random point in $ B^{\nu_\I}(x_0,\eps)$, which is a ball lying in the $n-1$-dimensional hyperplane  with the normal $\nu_\I$. Similarly if Player II wins the toss, then he tosses a biased coin. If he gets heads (with probability $\alpha$), the token is placed at $x_1=x_0+\nu_\II$. If he gets tails (with probability $\beta$), then the token is placed at a random point in $ B^{{\nu_\II}}(\eps,x_0)$. 
The game is played until the token exits $\Om$, and at the end Player II pays Player I the amount given by the pay-off function $F$. 
The \emph{value of the game}
is given by
\[
u(x_0)=\sup_{S_{\I}}\inf_{S_{\II}}\,\mathbb{E}_{S_{\I},S_{\II}}^{x_0}[F(x_\tau)].
\]

To avoid measurability problems, we can incorporate the continuous boundary correction in Section 4 of \cite{luirops14} as explained in \cite{hartikainen16}. However as we are studying local results, we only give a local form of the dynamic programming equation, and take it for granted.
By summing the possible outcomes of a single game round, we heuristically obtain  
\begin{equation}
\label{eq:dpp-tgw noise}
\begin{split}
u (x)=& \frac{1}{2} \sup_{ 0<\abs{\nu}\leq \eps}\bigg\{
\alpha u(x+\nu)  + \beta
\kint_{B^{\nu}_\eps} u (x+h) \ud \mathcal{L}^{n-1}(h)\bigg\}\\
& + \frac{1}{2}\inf_{0<\abs{\nu}\leq \eps}\bigg\{
\alpha u(x+\nu)  +\beta
\kint_{B^{\nu}_\eps} u (x+h) \ud \mathcal{L}^{n-1}(h)\bigg\},
\end{split}
\end{equation}
where $\mathcal L^{n-1}$ denotes the $n-1$-dimensional Lebesgue measure and $B^\nu_\eps:=B^\nu(0,\eps)$.
Similarly as in Section \ref{sec:p(x)}, we take the dynamic programming equation as our starting point straight away.

\subsection{Multidimensional dynamic programming}

We use the same $f$ as before, but now choose
\begin{equation}
\label{Caoletus}
C=\frac{10^{10}}{\delta^{2}\omega(\alpha)}\,.
\end{equation}

Fix $\eta>0$ and recall the counter assumption (\ref{vastaoletus}). Similarly as in Section \ref{sec:p=infty}, using
  (\ref{eq:dpp-tgw noise}), we end up with
\[
\begin{split}
M&\le u(x)-u(z)-f(x,z)+\eta
\\ &\le \sup_{\nu_x}
\frac{1}{2} \bigg\{
\alpha u(x+\nu_x)  + \beta
\kint_{B_\eps^{{\nu_x}}} u (x+h) \ud \mathcal{L}^{n-1}(h)\bigg\}\\
&\hspace{1 em} -\inf_{\nu_z}\frac{1}{2}\bigg\{
\alpha u(z+\nu_z)  +\beta
\kint_{B_\eps^{{\nu_z}}} u (z+h) \ud \mathcal{L}^{n-1}(h)\bigg\}\\
&\hspace{1 em} +
\inf_{\nu_x}\frac{1}{2}\bigg\{
\alpha u(x+\nu_x)  +\beta
\kint_{B_\eps^{{\nu_x}}} u (x+h) \ud \mathcal{L}^{n-1}(h)\bigg\}\\
&\hspace{1 em}-\sup_{\nu_z} \frac{1}{2} \bigg\{
\alpha u(z+\nu_z)  + \beta
\kint_{B_\eps^{{\nu_z}}} u (z+h) \ud \mathcal{L}^{n-1}(h)\bigg\}
-f(x,z)+\eta.
%
\end{split}
\]
Then, let  $P_{\nu_x,\nu_z}$ be a rotation, to be specified later,   from the hyperplane determined by $\nu_x$ to the hyperplane determined  by  $\nu_z$ so that we may briefly write
\[
\begin{split}
M\le
&\frac{1}{2} \sup_{\nu_x}\sup_{\nu_z}\bigg\{
\alpha( u(x+\nu_x)-u(z+\nu_z))  \\
&\hspace{1 em}+ \beta
\kint_{B_\eps^{{\nu_x}}} u (x+h)-u(z+P_{\nu_x,\nu_z}(h)) \ud \mathcal{L}^{n-1}(h)\bigg\}\\
&\hspace{1 em} + \frac{1}{2}\inf_{\nu_x}\inf_{\nu_z}\bigg\{
\alpha (u(x+\nu_x)-u(z+\nu_z))  \\
&\hspace{1 em}+\beta
\kint_{B_\eps^{{\nu_x}}} u (x+h)-u(z+ P_{\nu_x,\nu_z}(h))\ud \mathcal{L}^{n-1}(h)\bigg\}\\
&\hspace{1 em}-f(x,z)+\eta\\
&\le I_1+I_2-f(x,z)+\eta.
\end{split}
 \]
Next we choose $\nu^\I_x,\,\nu^\II_z\in \ol B(0,\eps)$ that give $\sup_{\nu_x}$ (denoted by $\nu^\I_x$) and $\sup_{\nu_z}$ in $I_1$ above up to an error $\eta/2>0$.  Then add and subtract $f$ similarly as in (\ref{eq:sup-inf}).
This gives
\[
\begin{split}
&I_1\le 
 \frac{M}{2} +\sup_{\nu_x,\nu_z}T(f,x,z,\nu_x,\nu_z)+\eta,
\end{split}
\]
where we used a shorthand notation
\[
\begin{split}
T f&:=T(f,x,z,\nu_x,\nu_z)\\
&:=\frac{\alpha}{2} f(x+\nu_x,z+\nu_z))+\frac{\beta}{2} \kint_{B_\eps^{{\nu_x}}} f(x+h,z+P_{\nu_x,\nu_z}(h)) \ud \mathcal{L}^{n-1}(h).
\end{split}
\]
Similarly to (\ref{eq:inf-sup}) we also  estimate
\[
\begin{split}
I_2
&\le \frac{M}{2}+\inf_{\nu_x,\nu_z}T(f,x,z,\nu_x,\nu_z)+\eta.
\end{split}
\]
In the previous two estimates, we used the counter assumption.

Similarly  as before, to obtain a contradiction,  it thus suffices to show 
\begin{equation}
\label{eq:f-p>1}
\begin{split}
\sup_{\nu_x,\nu_z}&T(f,x,z,\nu_x,\nu_z)
+\inf_{\nu_x,\nu_z}T(f,x,z,\nu_x,\nu_z)<f(x,z).
\end{split}
\end{equation}
We accomplish this in several steps.
\subsection{Estimates}
\subsubsection*{\bf Inequality \eqref{eq:f-p>1}, case $|x-z|> N\frac{\eps}{10}$}
In this case we may again focus our attention on $f_1$. We recall the Taylor expansion
\[
\begin{split}
f_1(x+h_x,z+h_z)=&f_1(x,z)+C\delta|x-z|^{\delta-1}(h_x-h_z)_V+2(x+z)\cdot (h_x+h_z)\\
&\,\,+\frac{C}{2}\delta|x-z|^{\delta-2}\big((\delta-1)(h_x-h_z)^2_V+(h_x-h_z)^2_{V^{\perp}}\big)\\
&\,\,+|h_x+h_z|^2\,+\mathcal{E}_{x,z}(h_x,h_z).
\end{split}
\]
Choose $  \nu'_x,  \nu'_z\in \ol B_\eps$ such that 
\[
\begin{split}
\sup_{\nu_x,\nu_z} T(f_1,x,z,\nu_x,\nu_z)\le T(f_1,x,z,  \nu'_x,  \nu'_z)+\eta.
\end{split}
\]

Let us first assume that
\begin{equation}\label{eq:vastaan}
(  \nu'_x-  \nu'_z)_V^2\geq (4-\theta)\eps^2\,\,.  
\end{equation}
Since $|  \nu'_x-  \nu'_z|\leq 2\eps$, this also means that 
\begin{equation}
\label{eq:vastaan-perp}
(  \nu'_x-  \nu'_z)^2_{V^{\perp}}\leq \theta\eps^2\,.
\end{equation}
In this case, the second order terms in the Taylor estimate of the $\alpha$-terms dominate  and thus yield the desired conclusion.
To be more precise, we estimate the quantity
\begin{equation}
\label{eq:main-estimate}
\begin{split}
T(&f_1,x,z,  \nu'_x,  \nu'_z)+T(f_1,x,z,-  \nu'_x,-  \nu'_z)-2f_1(x,z).
\end{split}
\end{equation}
To this end,
by using the above Taylor formula and the estimate (\ref{error2}) for the error term, we see that
the $\alpha$-terms in \eqref{eq:main-estimate} can be estimated as
\begin{equation}
\label{eq:alpha-bound}
\begin{split}
&\frac\alpha2  f_1(x+ \nu'_x,z+ \nu'_z))+\frac\alpha2 f_1(x- \nu'_x,z- \nu'_z))-2f_1(x,z)\\
&\le 0+\frac{C\alpha}{2}\delta|x-z|^{\delta-2}\big(\,2(\delta-1)(  \nu'_x-  \nu'_z)_{V}^2+2(  \nu'_x-  \nu'_z)^2_{V^{\perp}}   \,\big)\\
&\,\,\,\,\,\,+\alpha|  \nu'_x+  \nu'_z|^2+\frac\alpha2 \mathcal{E}_{x,z}(  \nu'_x,  \nu'_z)+\frac\alpha2 \mathcal{E}_{x,z}(-  \nu'_x,-  \nu'_z)\,\\
&\leq \frac{C\alpha }{2}\delta|x-z|^{\delta-2}\big(\,2(\delta-1)(4-\theta)\eps^2+2\theta\eps^2\,\big)\\
&\hspace{1 em}+\alpha (2\eps)^2\,+20\alpha \eps^2|x-z|^{\delta-2}.
\end{split}
\end{equation}
Next we estimate  the $\beta$-terms
\begin{equation}
\label{eq:beta-terms}
\begin{split}
\frac{\beta}{2}& \kint_{B_\eps^{{\nu'_x}}} f_1(x+h,z+P_{\nu'_x,\nu'_z}(h)) \ud \mathcal{L}^{n-1}(h)\\
&+\frac{\beta}{2} \kint_{B_\eps^{{-\nu'_x}}} f_1(x+h,z+P_{-\nu'_x,-\nu'_z}(h)) \ud \mathcal{L}^{n-1}(h)\\
\end{split}
\end{equation}
 in (\ref{eq:main-estimate}) by the Taylor series. By symmetry, the first order terms vanish. So far we have not used a particular form of the rotations $P_{\nu'_x,\nu'_z}$ and $P_{-\nu'_x,-\nu'_z}$. 
Now we choose  $P_{\nu'_x,\nu'_z}=P_{-\nu'_x,-\nu'_z}$ such that  $\abs{h-P_{\nu'_x,\nu'_z}(h)}^2\le \theta \eps^2 $. This and \eqref{eq:vastaan-perp} imply that  
\[
\begin{split}
C\delta|x-z|^{\delta-2}\kint_{B_\eps^{{  \nu_x}}}(h-P_{\nu'_x,\nu'_z}(h))^2_{V^{\perp}}\ud \mathcal L^{n-1}(h)\le C\delta|x-z|^{\delta-2} \theta \eps^2.
\end{split}
\] 
The same estimate naturally also holds for the second term of (\ref{eq:beta-terms}). We may also estimate 
$$
C(\delta-1)\delta|x-z|^{\delta-2}\kint_{B_\eps^{{  \nu'_x}}}(h-P_{\nu'_x,\nu'_z}(h))^2_{V}\ud \mathcal L^{n-1}(h)\le0
$$ 
and thus we may bound the quantity in (\ref{eq:beta-terms}) by
\begin{equation}
\label{eq:beta-bound}
\begin{split}
 C\beta \delta|x-z|^{\delta-2} \theta \eps^2+ \beta(2\eps)^2+\beta20 \eps^2|x-z|^{\delta-2}.
\end{split}
\end{equation}

Combining the estimates (\ref{eq:alpha-bound}) and (\ref{eq:beta-bound}), choosing small enough $\delta$ and $\theta$ as well as by using the assumption (\ref{Caoletus}) for $C$ along with choosing $\om(\alpha)$ small enough, we observe that   
\begin{equation}
 |x-z|^{\delta-2}(20+C\delta( \beta\theta-\alpha )) < -|x-z|^{\delta-2}10^8< -10^7\,.\notag
\end{equation}
This implies the negativity of the quantity in (\ref{eq:main-estimate}).

In turn, if (\ref{eq:vastaan}) above does not hold, implying that
\[
\begin{split}
(  \nu'_x-  \nu'_z)_V \leq (2-\frac{\theta}{4})\eps\,\,,  
\end{split}
\]
the desired estimate can be readily derived from the first order terms in the Taylor estimate. Indeed, the share of the second order  and error terms in the Taylor estimate of \eqref{eq:main-estimate} can be roughly estimated by (\ref{eq:error-and-second-rough}) as 
$$ \delta^2|x-z|^{\delta-1}\eps\,.$$
Furthermore, by using Taylor's formula and the above estimate for the second order terms and the error,  we have 
\begin{align*}
\sup_{}&Tf_1+\inf_{}Tf_1-2f_1(x,z)\\
\le&
T(f_1,x,z,  \nu'_x,  \nu'_z)+T(f_1,x,z,-\eps\frac{x-z}{|x-z|},\eps\frac{x-z}{|x-z|})-2f_1(x,z)\,+\eta\\
\leq &C\alpha\delta|x-z|^{\delta-1}((  \nu'_x-  \nu'_z)_V-2\eps)+4|x+z|2\eps+\delta^2|x-z|^{\delta-1}\eps\,+\eta \,\\
\leq & C\alpha \delta|x-z|^{\delta-1}(-\frac{\theta}{4}\eps)+16\eps+\delta^2|x-z|^{\delta-1}\eps\,+\eta\\
\leq &(\delta-\alpha \frac{\theta}{4} C)\delta|x-z|^{\delta-1}\eps+16\eps +\eta.
\end{align*}
Observe that above the first order terms in the Taylor series of $\beta$-terms vanish  in the integration by symmetry and also since $P_{\cdot,\cdot}$ is a rotation. 
By applying the assumption (\ref{Caoletus}) for $C$ for small enough $\om(\alpha)\le \alpha \theta$, it is easy to check the desired negativity for this quantity.  This completes the proof of the case $|x-z|>N\frac{\eps}{10}$. 


\subsubsection*{\bf Inequality (\ref{eq:f-p>1}), case $|x-z|\le  N\frac{\eps}{10}$}  
Recall that when $x,z\in B_1$, and $h_x,h_z\in \ol B_\eps$, then 
\[
\begin{split}
|f_1(x+h_x,z+h_z)-f_1(x,z)|\leq 3C\eps^{\delta}\,
\end{split}
\]
as in (\ref{perusperus}).
This implies
\[
\begin{split}
\sup Tf_1-f_1(x,z)\leq  3C\eps^{\delta}\,.
\end{split}
\]
Moreover, 
\[
\begin{split}
\sup T(f_1-f_2)\le \sup T(f_1-0)=\sup Tf_1.
\end{split}
\]
Letting $(i-1)\frac{\eps}{10}< |x-z|\le i\frac{\eps}{10}$, $i=1,2,\ldots,N$, we also have
\[
\begin{split}
\inf  T(f_1-f_2)&\le \sup Tf_1-\sup T f_2 \\
&\le \,\,\sup T f_1- \alpha C^{2(N-i+1)}\eps^{\delta}\\
&= \,\sup T f_1-\alpha (C^2-\frac2\alpha)C^{2(N-i)}\eps^{\delta}-2 C^{2(N-i)}\eps^{\delta}\,\\
&=\sup T  f_1 -\,\alpha  (C^2-\frac2\alpha)C^{2(N-i)}\eps^{\delta}-2 f_2(x,z)\,\\
&< \sup Tf_1-6C\eps^{\delta}-2 f_2(x,z)\,
\end{split}
\]
for small enough $\om(\alpha)$ in (\ref{Caoletus}).
Combining the previous estimates, 
we end up with 
\[
\begin{split}
\sup Tf+\inf T f&< \sup T f_1+\sup Tf_1 -6C\eps^\delta-2f_2(x,z)\\
&\le 2 f_1(x,z)+6C\eps^\delta  -6C\eps^\delta-2f_2(x,z)\\
&\le 2 f(x,z). 
\end{split}
\]
This completes the proof of Theorem \ref{thm:main}.




\def\cprime{$'$} \def\cprime{$'$}



\begin{thebibliography}{PSSW09}

\bibitem[AB10]{atarb10}
R.~Atar and A.~Budhiraja.
\newblock A stochastic differential game for the inhomogeneous $\infty$-{L}aplace
  equation.
\newblock {\em Ann. Probab.}, 38(2):498--531, 2010.


\bibitem[APSS12]{antunovicpss12}
T.~Antunovi{\'c}, Y.~Peres, S.~Sheffield, and S.~Somersille.
\newblock Tug-of-war and infinity {L}aplace equation with vanishing {N}eumann
  boundary condition.
\newblock {\em Comm. Partial Differential Equations}, 37(10):1839--1869, 2012.
%


\bibitem[BCF12a]{bjorklandcf12a}
C.~Bjorland, L.~Caffarelli, and A.~Figalli.
\newblock Non-local gradient dependent operators.
\newblock {\em Adv. Math.}, 230(4-6):1859--1894, 2012.

\bibitem[BCF12b]{bjorklandcf12}
C.~Bjorland, L.~Caffarelli, and A.~Figalli.
\newblock Nonlocal tug-of-war and the infinity fractional {L}aplacian.
\newblock {\em Comm. Pure Appl. Math.}, 65(3):337--380, 2012.

%
%
%
%
%


\bibitem[Har16]{hartikainen16}
H.~Hartikainen
\newblock  A dynamic programming principle with continuous solutions related to the p-Laplacian.
\newblock {\em  Differential Integral Equations}, 29(5--6):583--600, 2016.

\bibitem[IL90]{ishiil90}
H.~Ishii and P.-L. Lions.
\newblock Viscosity solutions of fully nonlinear second-order elliptic partial
  differential equations.
\newblock {\em J. Differential Equations}, 83(1):26--78, 1990.



\bibitem[JLP10]{juutinenlp10}
P.~Juutinen, T.~Lukkari, and M.~Parviainen.
\newblock Equivalence of viscosity and weak solutions for the $p(x)$-{L}aplacian.
\newblock 27(6):1471--1487, 2010.


\bibitem[Jul15]{julin15}
V.~Julin.
\newblock Generalized {H}arnack inequality for nonhomogeneous elliptic
  equations.
\newblock {\em Arch. Ration. Mech. Anal.}, 216(2):673--702, 2015.

\bibitem[KMP12]{kawohlmp12}
B.~Kawohl, J.~J. Manfredi, and M.~Parviainen.
\newblock Solutions of nonlinear {PDE}s in the sense of averages.
\newblock {\em J. Math. Pures Appl. (9)}, 97(2):173--188, 2012.
%
%
%

\bibitem[Kus15]{kusuoka15}
S.~Kusuoka.
\newblock H\"older continuity and bounds for fundamental solutions to
  nondivergence form parabolic equations.
\newblock {\em Anal. PDE}, 8(1):1--32, 2015.

\bibitem[LR86]{lindvallr93}
T.~Lindvall and L.~C.~G. Rogers.
\newblock Coupling of multidimensional diffusions by reflection.
\newblock {\em Ann. Probab.}, 14(3):860--872, 1986.


\bibitem[LPS13]{luirops13}
H.~Luiro, M.~Parviainen, and E.~Saksman.
\newblock Harnack's inequality for $p$-harmonic functions via stochastic games.
\newblock {\em Comm.\ Partial Differential Equations}, 38(11):1985--2003, 2013.

\bibitem[LPS14]{luirops14}
H.~Luiro, M.~Parviainen, and E.~Saksman.
\newblock On the existence and uniqueness of $p$-harmonious functions.
\newblock {\em Differential and Integral Equations}, 27(3/4):201--216, 2014.

\bibitem[LS15]{lius15}
Q.~Liu and A.~Schikorra.
\newblock General existence of solutions to dynamic programming equations.
\newblock {\em Commun. Pure Appl. Anal.}, 14(1):167--184, 2015.

\bibitem[MPR10]{manfredipr10}
J.J. Manfredi, M.~Parviainen, and J.D. Rossi.
\newblock An asymptotic mean value characterization for $p$-harmonic functions.
\newblock {\em Proc. Amer. Math. Soc.}, 258:713--728, 2010.

\bibitem[MPR12]{manfredipr12}
J.J. Manfredi, M.~Parviainen, and J.D. Rossi.
\newblock On the definition and properties of p-harmonious functions.
\newblock {\em Ann. Scuola Norm. Sup. Pisa Cl. Sci.}, 11(2):215--241, 2012.

%
\bibitem[Obe05]{oberman05}
A.~M. Oberman.
\newblock A convergent difference scheme for the infinity {L}aplacian:
  construction of absolutely minimizing {L}ipschitz extensions.
\newblock {\em Math. Comp.}, 74(251):1217--1230, 2005.

\bibitem[PS08]{peress08}
Y.~Peres and S.~Sheffield.
\newblock Tug-of-war with noise: a game-theoretic view of the
  {$p$}-{L}aplacian.
\newblock {\em Duke Math. J.}, 145(1):91--120, 2008.

\bibitem[PSSW09]{peresssw09}
Y.~Peres, O.~Schramm, S.~Sheffield, and D.~B. Wilson.
\newblock Tug-of-war and the infinity {L}aplacian.
\newblock {\em J. Amer. Math. Soc.}, 22(1):167--210, 2009.

\bibitem[PP13]{porrettap13}
A.~Porretta and E.~Priola.
\newblock Global {L}ipschitz regularizing effects for linear and nonlinear
  parabolic equations.
\newblock {\em J. Math. Pures Appl. (9)}, 100(5):633--686, 2013.


\bibitem[PW06]{priolaw06}
E.~Priola and F.-Y. Wang.
\newblock Gradient estimates for diffusion semigroups with singular
  coefficients.
\newblock {\em J. Funct. Anal.}, 236(1):244--264, 2006.

\bibitem[Ruo16]{ruosteenoja16}
E.~Ruosteenoja.
\newblock Local regularity results for value functions of tug-of-war with noise
  and running payoff.
\newblock {\em Adv. Calc. Var.}, 9(1):1--17, 2016.

\end{thebibliography}

\end{document}